\documentclass{amsart}
\usepackage{amsfonts,amssymb}
\usepackage{graphicx}
\usepackage[noadjust]{cite}

\usepackage{amsmath,bbm}

\usepackage[hidelinks]{hyperref}
\usepackage[T1]{fontenc}
\usepackage[cp1250]{inputenc}
\usepackage{float}
\usepackage{enumitem}
\usepackage{mathabx}

\usepackage{upgreek}

\newcommand{\assign}{:=}

\newcommand{\tmem}[1]{{\em #1\/}}

\newcommand{\tmop}[1]{\ensuremath{\operatorname{#1}}}

\newtheorem{lem}{Lemma}[section]
\newtheorem{cor}[lem]{Corollary}

\newtheorem{thm}[lem]{Theorem}
\newtheorem{prop}[lem]{Proposition}

\newtheorem{rk}{Remark}
\newtheorem*{problem**}{Problem}
\newtheorem*{thm**}{Theorem}
\newtheorem{ex}{Example}
\newtheorem*{ex**}{Example}
\newtheorem*{rk**}{Remark}

\newenvironment{prf}{\noindent\textbf{Proof.\ }}{\hspace*{\fill}$\Box$\medskip}

\title[Jumps of Milnor numbers  of Brieskorn-Pham singularities]{Jumps of Milnor numbers  of Brieskorn-Pham singularities in non-degenerate families}
\author[T. Krasi{\'n}ski and J. Walewska]{Tadeusz Krasi{\'n}ski and Justyna Walewska}
\address{Faculty of Mathematics and Computer Science, University of \L\'od\'z, ul. Banacha 22, 90--238 \L\'od\'z, Poland}
\email[Tadeusz Krasi\'nski]{krasinsk@uni.lodz.pl}
\email[Justyna Walewska]{walewska@math.uni.lodz.pl}

\begin{document}

	\renewcommand{\thepage}{[\arabic{page}]}

\begin{abstract}
The jump of the Milnor number of an isolated singularity $f_{0}$ is the
minimal non-zero difference between the Milnor numbers of $f_{0}$ and one of
its deformation $(f_{s}).$ In the case $f_{s}$ are non-degenerate
singularities we call the jump non-degenerate. We give a formula (an
inductive algorithm using diophantine equations) for the non-degenerate jump of $f_{0}$ in the case
$f_{0}$ is a convenient singularity with only one $(n-1)$-dimensional face of
its Newton diagram which equivalently (in our problem) can be replaced by the Brieskorn-Pham singularities.

\end{abstract}

\maketitle

\section{Introduction}\label{Par 1}

Let $f_0 : (\mathbbm{C}^n, 0) \rightarrow (\mathbbm{C}, 0)$ be an
{\tmem{(isolated) singularity}}, i.e. $f_0$ is the germ at $0$ of a
holomorphic function having an isolated critical point at $0 \in
\mathbbm{C}^n$, and $0 \in \mathbbm{C}$ as the corresponding critical value.
More specifically, there exists a representative $\hat{f}_0 : U \rightarrow
\mathbbm{C}$ of $f_0$ holomorphic in an open neighborhood $U$ of the point $0
\in \mathbbm{C}^n$ such that:

\begin{itemize}
	\item[$\bullet$] $\hat{f}_0 (0) = 0$,
	
	\item[$\bullet$] $\nabla \hat{f}_0 (0) = 0$,
	
	\item[$\bullet$] $\nabla \hat{f}_0 (z) \neq 0$ for $z \in U\backslash \{ 0 \}$,
\end{itemize}

{\noindent}where for a holomorphic function $f$ we put $\nabla f \assign
(\partial f / \partial z_1, \ldots, \partial f / \partial z_n)$.

The ring of germs of holomorphic functions of $n$ variables will be denoted by $\mathcal{O}_n$.

A {\tmem{deformation of the singularity}} $f_0$ is the germ of a holomorphic
function $f = f (s, z) : (\mathbbm{C} \times \mathbbm{C}^n, 0) \rightarrow
(\mathbbm{C}, 0)$ such that:

\begin{itemize}
	\item[$\bullet$] $f (0, z) = f_0(z)$,
	
	\item[$\bullet$] $f (s, 0) = 0$,

\end{itemize}
The deformation $f (s, z)$ of the singularity $f_0$ will also be treated as a
family $(f_s)$ of germs, putting $f_s (z) \assign f (s, z)$. Since $f_0$ is an isolated singularity, $f_s$ has also isolated singularities near the origin, for sufficiently small $s$ {\cite[Theorem 2.6 in Chap.~I]{GLS07}}.
\begin{rk}
	\label{nota1}Notice that in the deformation $(f_s)$ there can occur in
	particular smooth germs, that is germs satisfying $\nabla f_s (0)
	\neq 0$. In this context,
	the symbol $\nabla f_s$ will always denote $\nabla_z f_s(z)$.
	
\end{rk}

By the above assumptions, for every sufficiently small $s$,
one can define a (finite) number ${\upmu}_s$ as the Milnor number of $f_s$,
namely
\[ {\upmu}_s \assign {\upmu} (f_s) \mathrel{=} \dim_{\mathbbm{C}} 
\mathbin{\mathbin{\mathcal{O}}_n / \mathbin{(\nabla f_s)}}
\mathord{\mathrel{}} \mathrel \text{,}
\]
where $(\nabla f_s)$ is the ideal in $\mathcal{O}_n$ generated by $\partial f_s/\partial z_1, \ldots,\partial f_s/\partial z_n$.

Since the Milnor number is upper semi-continuous in the Zariski topology \\
{\cite[Prop. 2.57 in Chap.~II]{GLS07}}, there exists an open neighborhood
$S$ of the point $0 \in \mathbbm{C}$ such that

\begin{itemize}
	\item[$\bullet$] ${\upmu}_s = \tmop{const.}$ for $s \in S \setminus \{ 0 \}$,
	
	\item[$\bullet$] ${\upmu}_0 \geqslant {\upmu}_s$ for $s \in S$.
\end{itemize}

The (constant) difference ${\upmu}_0 -{\upmu}_s$ for $s \in S \setminus \{
0 \}$ will be called {\tmem{the jump of the deformation $(f_s)$}} and denoted
by $\uplambda ((f_s))$. The smallest non-zero value among the jumps of all
deformations of the singularity $f_0$ will be called {\tmem{the jump
		(of the Milnor number) of the singularity $f_0$}} and denoted by $\uplambda
(f_0)$.

The first general result concerning the jump was S.~Guse{\u\i}n-Zade's
{\cite{Gus93}}, who proved that there exist singularities $f_0$ for which
$\uplambda (f_0) > 1$ and that for irreducible plane curve singularities it
holds $\uplambda (f_0) = 1$. In {\cite{BK14}} the authors proved that $\uplambda
(f_0)$ is not a topological invariant of $f_0$ but it is an invariant of the
stable equivalence of singularities. The computation of $\uplambda (f_0)$ is not an
easy task. It is related to the problem of adjacency of classes of singularities. Only for few classes of singularities we know the exact value of
$\uplambda (f_0)$. For plane curve singularities ($n=2$) in the $X_9$ class (see {\cite{AGZV85}} for terminology) that is singularities of the form $f^{a}_0(x,y):=x^4+y^4+ax^2y^2$, $a\in\mathbb{C}$, $a^2\not=4$, we have $\lambda(f^{a}_0)=2$ (see {\cite{BK14}}) and for specific homogeneous singularities $f^{d}_0(x,y)=x^d+y^d$, $d\geqslant 2$, we have $\lambda(f^{d}_0)=\left[d/2\right]$, where
$[x]$ is the integer part of $x$ (see {\cite{BKW14}}).

In the paper we consider a weaker problem: {\tmem{compute the jump
		$\uplambda^{\tmop{nd}} (f_0)$ of $f_0$ over all non-degenerate deformations of
		$f_0$}} (i.e.~the $f_s$ in deformations $(f_s)$ of $f_0$ are
non-degenerate singularities). Clearly, we always have $\uplambda (f_0)
\leqslant \uplambda^{\tmop{nd}} (f_0)$. Up to now, this problem has been studied
only for plane curve singularities by A.~Bodin {\cite{Bod07}}, S.~Brzostowski, T.~Krasi\'nski and J.~Walewska {\cite{Wal10}}, \cite{Wal13}, {\cite{BKW14}}. 
We give a formula (more precisely an algorithm) in $n$-dimensional case for $\lambda^{\operatorname{nd}}(f_0)$, where $f_0$ is non-degenerate, convenient and has its Newton diagram reduced to one $(n-1)$-dimensional face. Since for non-degenerate singularities $f_s$ of a deformation $(f_s)$ of $f_0$ the Milnor number $\upmu(f_s)$ depends only on the Newton diagram of $f_s$ and not on specific coefficients of $f_s$ (see Section 2) we may restrict considerations to the class of the Brieskorn-Pham singularities
\begin{equation}\label{rownanie osobliwosci Brieskorna-Phama}
f_0(\textbf{z})=z_1^{p_1}+\ldots+z_n^{p_n},\quad p_i\geqslant 2, \quad i=1,\ldots,n,
\end{equation}
where $\textbf{z}=(z_1,\ldots,z_n)$.

\begin{rk}
	Let us note that we compute non-degenerate jump of a singularity in a fixed system of coordinates. The reason is that non-degeneracy depends on coordinate system. For instance, for singularity $f_0(x,y)=x^3+y^3$, $\lambda^{\operatorname{nd}}(f_0)=2$ in the original system of coordinates while for $L: x=x'-y'$, $y=y'$, we have $\lambda^{\operatorname{nd}}(f_0\circ L)=1$.
	
	The problem of finding a minimal non-degenerate jump of a singularity in arbitrary system of coordinates is open.
\end{rk}

\section{Non-degenerate singularities}

In this Section we recall the notion of non-degenerate singularities. Let
$f_0(\textbf{z})=\sum_{\textbf{i}\in\mathbb{N}_0^{n}}a_{\textbf{i}}\textbf{z}^{\textbf{i}}$
be a singularity and
$\tmop{supp}(f_0)=\{\textbf{i}\in\mathbb{N}_0^{n}:a_{\textbf{i}}\not=0\}$
the {\tmem{support of $f_0$}}. The {\tmem{Newton polyhedron $\Gamma_+ (f_0)\subset \mathbb{R}^n$}} of $f_0$ is the convex hull of the set 
$ \bigcup_{\textbf{i} \in \tmop{supp} (f_0)} (\textbf{i} +\mathbbm{R}_+^n). $
The finite set of compact faces (of all dimensions) of the boundary of $\Gamma_+ (f_0)$ constitutes the \tmem{Newton diagram} of $f_0$ and is denoted by $\Gamma (f_0)$. The singularity $f_0$ is \tmem{convenient} if $\Gamma_+ (f_0)$ intersects all coordinate axes in $\mathbb{R}^{n}$. For each face $S\in \Gamma(f_0)$ we define a weighted homogeneous polynomial
$(f_0)_S(\textbf{z})=\sum_{\textbf{i}\in S}a_{\textbf{i}}\textbf{z}^{\textbf{i}}.$
We call the singularity $f_0$ \tmem{non-degenerate on} $S\in \Gamma(f_0)$ if the system of equations
$\partial(f_0)_S/\partial z_i(\textbf{z})=0$, $i=1,\ldots,n$
has no solutions in $(\mathbbm{C}^{*})^n$. The singularity $f_0$ is \tmem{non-degenerate} if $f_0$ is non-degenerate on every face $S\in \Gamma(f_0)$.

In the sequel we put the following notations for coordinate hyperplanes in $\mathbb{R}^n$. Let $I=\{i_1,\ldots,i_k\}\subset\{1,2,\ldots,n\}$. We denote $$Ox_{i_1}\ldots x_{i_k}:=\{\mathbbm{x}=(x_1,\ldots,x_n):x_j=0,\; j\not\in I\}$$ and $$\{x_{i_1}=\ldots=x_{i_k}=0\}:=\{\mathbbm{x}=(x_1,\ldots,x_n):x_{i_1},\ldots=x_{i_k}=0\}.$$ For a convenient singularity $f_0$ we define the \tmem{Newton number} $\upnu(f_0)$ of $f_0$ by
$$\upnu(f_0)=n!V-\sum_{i=1}^{n}(n-1)!V_i+\sum_{\substack{
		i,j=1\\
		i<j}
}^{n}(n-2)!V_{ij}-\ldots+(-1)^{n},$$
where $V$ is the $n$-dimensional volume under $\Gamma_{+}(f_0)$, $V_i$ is the $(n-1)$-dimensional volume under $\Gamma_{+}(f_0)$ on the hyperplane $\{x_i=0\}$, $V_{ij}$ is the $(n-2)$-dimensional volume under $\Gamma_{+}(f_0)$ on the hyperplane $\{x_i=x_j=0\}$ and so on.

The importance of $\upnu (f_0)$ has its source in the celebrated {\it{Kouchnirenko
		theorem}}.

\begin{thm}[\cite{Kou76}]If $f_0$ is a convenient
	singularity, then
	\begin{enumerate}
		\item\label{warK1} ${\upmu} (f_0) \geqslant \upnu (f_0)$,
		
		\item\label{warK2} if $f_0$ is non-degenerate then ${\upmu} (f_0) = \upnu (f_0)$.
	\end{enumerate}
\end{thm}
\begin{rk}
	For non-convenient singularities the Kouchnirenko Theorem also holds provided $\nu(f_0)$ is appropriately defined (see for instance \textup{\cite{BO16}}). It is interesting that both conditions of item \ref{warK2}. of the Theorem are equivalent in $2$-dimensional case by A.~P\l oski \textup{\cite{Plo99}}. For weaker non-degeneracy conditions equivalent to the equality $\mu(f_0)=\nu(f_0)$ in $n$-dimensional case see recent preprint by P.~Mondal \textup{\cite{Mon16}}.
\end{rk}

\section{Non-degenerate jump of Milnor numbers of singularities}
Let $f_0\in\mathcal{O}_n$ be a singularity. A deformation $(f_s)$ of $f_0$ is called \textit{non-degenerate} if $f_s$ is non-degenerate for  $s\neq 0$. The set of all non-degenerate deformations of the singularity $f_0$ will be denoted by $\mathcal{D}^{\operatorname{nd}}(f_0)$. \textit{Non-degenerate jump} $\uplambda^{\operatorname{nd}}(f_0)$ \textit{of the singularity} $f_0$ is the minimal of non-zero jumps over all non-degenerate deformations of $f_0$, which means
$$\uplambda^{\operatorname{nd}}(f_0):=\min_{(f_s)\in \mathcal{D}^{\operatorname{nd}}_0 (f_0)} \uplambda ((f_s)),$$
where by $\mathcal{D}^{\operatorname{nd}}_0 (f_0)$ we denote all non-degenerate deformations $(f_s)$ of $f_0$ for which $\uplambda ((f_s))\neq 0$.

Obviously
\begin{prop}
	For each singularity $f_0$ we have the inequality 
	$$\uplambda (f_0)\leq \uplambda^{\operatorname{nd}}(f_0).$$
\end{prop}

From the Kouchnirenko Theorem we get

\begin{cor}\label{wnio1}
	If $f_0$ and $\tilde{f_0}$ are non-degenerate and convenient singularities and \\$\Gamma(f_0)=\Gamma(\tilde{f_0})$ then $\uplambda^{\operatorname{nd}}(f_0)=\uplambda^{\operatorname{nd}}(\tilde{f_0})$.
\end{cor}

\begin{cor}\label{wnio2}
	If $f_0$ is a convenient and non-degenerate singularity such that the Newton diagram $\Gamma(f_0)$ has only one $(n-1)$-dimensional face which intersects the axes $Ox_1,\ldots,Ox_n$ in points $(p_1,0,\ldots,0),\ldots, (0,\ldots,0,p_n)$ respectively, then for the Brieskorn-Pham singularity $\bar{f}(\textbf{z})=z_1^{p_1}+\ldots+z_n^{p_n}$ we have 
	$$\lambda^{\operatorname{nd}}(f_0)=\lambda^{\operatorname{nd}}(\bar{f}).$$
\end{cor}

In investigations concerning $\uplambda^{\operatorname{nd}}(f_0)$ we may restrict our attention to non-dege-\ nerate $f_0$ because the non-degenerate jump for degenerate singularities can be found using the proposition below (cf.~\cite[Lemma 5]{Bod07}). For a given $f_0$ let $f^{\operatorname{nd}}_0$ denote any non-degenerate singularity for which $\Gamma(f_0)=\Gamma(f^{\operatorname{nd}}_0)$. Such singularities always exist.

\begin{prop}\label{wzorklam}
	If $f_0$ is degenerate then
	
	$$\uplambda^{\operatorname{nd}}(f_0)=\bigg\{
	\begin{array}{ll}
	\upmu(f_0)-\upmu(f^{\operatorname{nd}}_0), & \text{if }\upmu(f_0)-\upmu(f^{\operatorname{nd}}_0)>0\\
	\uplambda^{\operatorname{nd}}(f^{\operatorname{nd}}_0), & \text{if }\upmu(f_0)-\upmu(f^{\operatorname{nd}}_0)=0
	\end{array}.$$
\end{prop}
\begin{prf}
	This follows from the fact that a generic small perturbation of coefficients of these monomials of $f_0$ which correspond to points belonging to $\bigcup\Gamma(f_0)$ (which are finite in number) give us non-degenerate singularities with the same Newton polyhedron as $f_0$.
\end{prf}

\begin{rk}
	By the P\l oski theorem \textup{(\cite[Lemma 2.2]{Plo90}, \cite[Theorem 1.1]{Plo99})}, for degenerate plane curve singularities ($n=2$) the second possibility in Proposition \ref{wzorklam} is excluded.
\end{rk}

A crucial role in the search for the formula for $\uplambda^{\operatorname{nd}}(f_0)$ will be played by the monotonicity of the Newton number with respect to the Newton polyhedron. Namely, J.~Gwo\'zdziewicz \cite{Gwo08}, M.~Furuya \cite{Fur04} and C.~Bivi\'a-Ausina \cite{B-A09} proved.

\begin{thm}[Monotonicity Theorem]
	Let $f_0, \tilde{f_0}\in\mathcal{O}_n$ be two convenient singularities such that $\Gamma_{+}(f_0)\subset \Gamma_{+}(\tilde{f_0})$. Then $\upnu(f_0)\geqslant \upnu(\tilde{f_0})$.
\end{thm}

 We need a stronger version of the above theorem in the particular case $f_0$ has only one $(n-1)$-dimensional face and $\tilde{f_0}$ differs from $f_0$ in one monomial. 

\begin{thm}\label{stronger monotonicity theorem}
	If $f_0$ is a convenient and non-degenerate singularity such that the Newton diagram $\Gamma(f_0)$ has only one $(n-1)$-dimensional face and $f^{\textup{\textbf{i}}}_0(\textbf{z})\assign f_0(\textbf{z})+\textbf{z}^{\textup{\textup{\textbf{i}}}}$, where the lattice point  $\textup{\textbf{i}}\not=0$ lies under $\Gamma_+ (f_0)$ i.e. $0\not=\textup{\textbf{i}}\in\mathbb{N}_0^{n}\setminus\Gamma_+(f_0)$, then $\upnu(f_0)>\upnu(f^{\textup{\textbf{i}}}_0)$.
\end{thm}

\begin{prf}
	By the Kouchnirenko Theorem we may consider only the Brieskorn-Pham singularities (\ref{rownanie osobliwosci Brieskorna-Phama}).	
	We will use induction with respect to the number $n$ of variables. For $n=1$ we have $f_0(\textbf{z})=z^{p_1}$, $p_1\geqslant 2$, and $f^{\textbf{i}}_0(\textbf{z})=z^{p_1}+z^i$, $1\leqslant i<p_1$. Hence $\upnu(f_0)=p_1-1$ and $\upnu(f^{\textbf{i}}_0)=i-1$. This gives $\upnu(f_0)>\upnu(f^{\textbf{i}}_0)$. 
	
	Suppose the theorem is true for any Brieskorn-Pham singularity in $(n-1)$ variables. Let now $f_0(\textbf{z})=z_1^{p_1}+\ldots+z_n^{p_n}$, $p_i\geqslant 2$, $i=1,\ldots,n$. Let $P_i=(0,\ldots,p_i,\ldots,0)$.
	
	If $\textbf{i}$ does not lie on any coordinate hyperplane $\{x_i=0\}$ then $\Gamma_+ (f_0)$ and $\Gamma_+ (f^{\textbf{i}}_0)$ are identical on all these hyperplanes (see Figure \ref{rys:1}(a)) and hence
	$$\upnu(f_0)-\upnu(f^{\textbf{i}}_0)=n!\operatorname{vol}(\Delta(\textbf{i},P_1,\ldots,P_n))>0,$$
	where $\Delta(\textbf{i},P_1,\ldots,P_n)$ is the $n$-dimensional simplex with vertices $\textbf{i},P_1,\ldots,P_n$.

\begin{figure}[H]
	\centering
	\includegraphics[width=\textwidth]{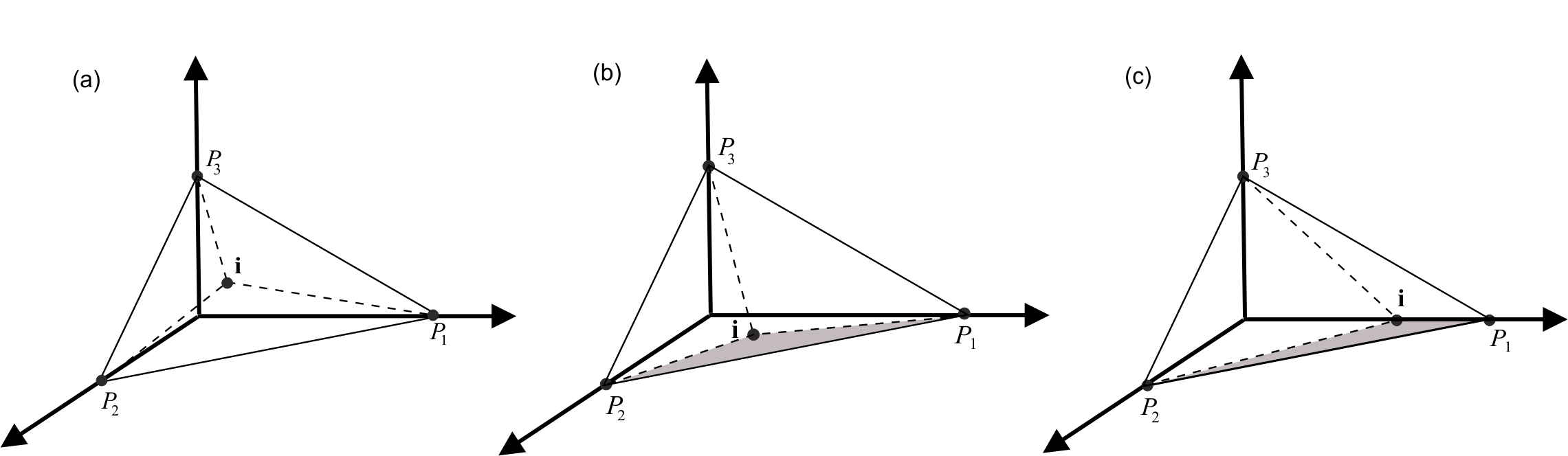}\\
	\caption{Possible cases in dimension three.}\label{rys:1}
\end{figure}

Take now a lattice point $\textbf{i}\not=0$ under $\Gamma_+ (f_0)$ which lies on a coordinate hyperplane. For simplicity we may assume $\textbf{i}\in Ox_1\ldots x_{n-1}$ i.e. $\textbf{i}=(i_1,\ldots,i_{n-1},0)$. Two possible cases in dimension three are illustrated in Figure \ref{rys:1}(b),(c). 
By definition of the Newton number the difference $\upnu(f_0)-\upnu(f^{\textbf{i}}_0)$ is equal to the sum of volumes of differences of $\Gamma_+ (f_0)$ and $\Gamma_+ (f^{\textbf{i}}_0)$ on all hyperplanes $Ox_{i_1}\ldots x_{i_l}$ multiplied by $(-1)^{n-l}l!$ $({i_1},\ldots ,{i_l} \text{ are different and } l=1,\ldots,n)$. First we consider these differencies in hyperplanes which do not contain the axis $Ox_n$. They are equal to the sum of volumes of differencies for truncated singularities $^nf_0(\textbf{z})\assign z_1^{p_1}+\ldots+z_{n-1}^{p_{n-1}}$ and $^nf^{\textbf{i}}_0(\textbf{z})\assign    {^nf_0(\textbf{z})}+\textbf{z}^{\textbf{i}}$ but taken with opposite signs. So, the sum of this differencies is equal to $-(\upnu(^nf_0)-\upnu(^nf^{\textbf{i}}_0))$. Now we consider differencies in coordinate hyperplanes which contain the axis $Ox_n$. Each such difference in a coordinate hyperplane $Ox_{i_1}\ldots x_{i_k}x_n$ is a pyramid whose base is equal to the difference, say $W$, of $\Gamma_+ (f_0)$ and $\Gamma_+ (f^{\textbf{i}}_0)$ in $Ox_{i_1}\ldots x_{i_k}$ and the vertex is equal to $P_n=(0,\ldots,0,p_n)$. Hence, in $\upnu(f_0)-\upnu(f^{\textbf{i}}_0)$ this pyramid will give the term 
	$$(-1)^{n-(k+1)}(k+1)!\dfrac{\operatorname{vol}(W)p_n}{k+1}=(-1)^{n-1-k}k!\operatorname{vol}(W)p_n.$$
	But the latter expression is precisely the term in $(\upnu(^nf_0)-\upnu(^nf^{\textbf{i}}_0))p_n$ associated to the difference $W$. Together we get $$\upnu(f_0)-\upnu(f^{\textbf{i}}_0)=(\upnu(^nf_0)-\upnu(^nf^{\textbf{i}}_0))(p_n-1).$$
	Since $p_n>1$ and by induction hypothesis $\upnu(^nf_0)-\upnu(^nf^{\textbf{i}}_0)$ is positive, we get the assertion of the theorem.
\end{prf}

Now we define specific non-degenerate deformations of a convenient and non-degenerate singularity $f_0\in\mathcal{O}_n$ (in particular for Brieskorn-Pham singularities). Denote by $J$ the set of lattice points (points with integer coordinates) $\textbf{i}\not=0$ lying under the Newton diagram of $f_0$ i.e. $0\not=\textbf{i}\in\mathbb{N}_0^n\setminus\Gamma_+(f_0)$. Obviously, $J$ is a finite set. For $\textbf{i}=(i_1,\ldots,i_n)\in J$ we define the \tmem{monomial deformation} $(f_s^{\textbf{i}})_{s\in\mathbb{C}}$ of $f_0$ by the formula
$$f^{\textbf{i}}_s(\textbf{z}):=f_0(\textbf{z})+s\textbf{z}^{\textbf{i}}.$$

\begin{prop}
	For every $\textup{\textbf{i}}\in J$ the deformation $(f^{\textup{\textbf{i}}}_s)$ of $f_0$ is convenient and non-degenerate for all sufficiently small $|s|$.
\end{prop}

\begin{prf}
	See \cite{Kou76} or \cite[Appendix]{Oka79}.
	
\end{prf}

Combining the Monotonicity Theorem with the above proposition we reach the conclusion that in order to find $\uplambda^{\operatorname{nd}}(f_0)$ it is enough to consider only the non-degenerate deformations of the type $(f^{\textup{\textbf{i}}}_s)$.
\begin{thm}\label{twieA}
	If $f_0$ is a convenient and non-degenerate singularity such that the Newton diagram $\Gamma(f_0)$ has only one $(n-1)$-dimensional face, then
	$$\uplambda^{\operatorname{nd}}(f_0)=\min_{\textup{\textbf{i}} \in J}\uplambda((f^{\textup{\textbf{i}}}_s)).$$
\end{thm}

\begin{prf}
	By the Kouchnirenko theorem it suffices to consider non-degenerate deformations of $f_0$ of the form
	
	\begin{equation}\label{postac deformacji z n punktami}
	f_s(\textbf{z})=f_0(\textbf{z}) + \sum_{\textup{\textbf{i}}\in J}a_{\textup{\textbf{i}}}(s)z^{\textup{\textbf{i}}},
	\end{equation}
	where $a_{\textup{\textbf{i}}}(s)$ are holomorphic at $0\in\mathbbm{C}$ and $a_{\textup{\textbf{i}}}(0)=0$. Then by the Monotonicity Theorem and Theorem \ref{stronger monotonicity theorem} we may restrict the scope of deformations (\ref{postac deformacji z n punktami}) to deformations with only one term added i.e.~to the deformations $(f^{\textup{\textbf{i}}}_s)$ for $\textup{\textbf{i}}\in J$.
\end{prf}

\section{An algorithm for the non-degenerate jump of Brieskorn-Pham singularities} 

In this Section we give an algorithm for calculating $\uplambda^{\operatorname{nd}}(f_0)$ for a non-degenerate and convenient singularity with one $(n-1)$-dimensional face of its Newton diagram.  By Corollary \ref{wnio2} we may assume that $f_0$ is a Brieskorn-Pham singularity
\begin{equation}
f_0(\textbf{z})=z_1^{p_1}+\ldots+z_n^{p_n}, p_i\geqslant 2.
\end{equation}

We have $\upmu(f_0)=\prod\limits_{i=1}^{n}(p_i-1)$. By Theorem \ref{twieA} to find $\lambda^{\operatorname{nd}}(f_0)$ we have to compare the jumps of all monomial deformations $(f^{\textbf{i}}_s)_{s\in\mathbb{C}}$ where $\textbf{i}\not=0$ is a lattice point lying under $\Gamma_+(f_0)$. First we prove a lemma relating the jumps of monomial deformations of the Brieskorn-Pham singularities to the jumps of their "truncations".

\begin{lem}\label{wazny lemat}
	Let $f_0(\textbf{z})=z_1^{p_1}+\ldots+z_n^{p_n}$ be a Brieskorn-Pham singularity. Define $^kf_0(\textbf{z})\assign z_1^{p_1}+\ldots+z_{k-1}^{p_{k-1}}+z_{k+1}^{p_{k+1}}+\ldots+z_n^{p_n}$, $k=1,\ldots,n$. Let $\textup{\textbf{i}}\not=0$ be a lattice point under $\Gamma_+(f_0)$ lying in the hyperplane $\{x_k=0\}$. Then
	
	\begin{equation}\label{rownanie w waznym lemacie}
	\lambda((f_s^{\textup{\textbf{i}}}))=\lambda((^kf_s^{\textup{\textbf{i}}}))(p_k-1).
	\end{equation}
\end{lem}

\begin{prf}
	The proof follows the line of the proof of Theorem \ref{stronger monotonicity theorem}. For simplicity of notations we take $k=n$, i.e. $\textbf{i}=(i_1,\ldots,i_{n-1},0)$. Since $\lambda((f_s^{\textbf{i}}))=\upnu(f_0)-\upnu(f_s^{\textbf{i}})$, $s\not=0$, by definition of the Newton number the difference $\upnu(f_0)-\upnu(f^{\textbf{i}}_s)$ is equal to the sum of volumes of differences of $\Gamma_+ (f_0)$ and $\Gamma_+ (f^{\textbf{i}}_s)$ on all hyperplanes $Ox_{i_1}\ldots x_{i_l}$ multiplied by $(-1)^{n-l}l!$ $({i_1},\ldots,{i_l} \text{ are different and } l=1,\ldots,n)$. 
	Similarly as in the proof of Theorem \ref{stronger monotonicity theorem} we divide all differencies in hyperplanes into two kinds:
	in those hyperplanes which do not contain the axis $Ox_n$ and those which do contain the axis $Ox_n$. The sum of the first ones is equal to the number $-\lambda((^nf_s^{\textbf{i}}))$ while the second ones $\lambda((^nf_s^{\textbf{i}}))p_n$. Together we get equality (\ref{rownanie w waznym lemacie}) for $k=n$. 
\end{prf}

Now we present an inductive algorithm for $\lambda^{\operatorname{nd}}(f_0)$ where $f_0$ is any Brieskorn-Pham singularity
$$f_0(\textbf{z})=z_1^{p_1}+\ldots+z_n^{p_n},\quad p_i\geqslant 2,\quad i=1,\ldots,n$$
with respect to the number of variables $n$. In each step we give the formula for $\lambda^{\operatorname{nd}}(f_0)$ and simultaneously the point $\textbf{i}\in J$ (precisely the deformation $(f^{\textbf{i}}_s)$) realizing this jump. For $n=1$ we have $f_0(z_1)=z_1^{p_1}$, $p_1\geqslant 2$ and $J=\{1,2,\ldots,p_1-1\}$. We have $\lambda^{\operatorname{nd}}(f_0)=1$ and the point $\textbf{i}=(p_1-1)$ realizes this jump.

Assume that for all Brieskorn-Pham singularities $f_0$ in $(n-1)$-variables we can compute $\lambda^{\operatorname{nd}}(f_0)$ and can indicate the point $\textbf{i}$ realizing this jump.

Take the Brieskorn-Pham singularity $f_0(\textbf{z})=z_1^{p_1}+\ldots+z_n^{p_n}, p_i\geqslant 2, n\geqslant 2$, in $n$ variables. By Theorem \ref{twieA} we have to compare the jumps of all monomial deformations $(f_s^{\textbf{i}})$, where $\textbf{i}\not=0$ is a lattice point lying under $\Gamma_+(f_0)$.

\medskip
A. First we consider points $\textbf{i}\not=0$ lying in coordinate hyperplanes of dimension $(n-1)$. By induction hypotheses on each coordinate hyperplane \\$Ox_1\ldots x_{k-1}x_{k+1}\ldots x_n$ there exists a point $\textbf{i}_k$ realizing non-degenerate jump of the truncated singularity $^kf_0=z_1^{p_1}+\ldots+z_{k-1}^{p_{k-1}}+z_{k+1}^{p_{k+1}}+\ldots+z_n^{p_n}$ in $(n-1)$ variables. Denote this jump by $\lambda_k\assign\lambda^{\operatorname{nd}}(^kf_0)=\lambda((^kf_s^{\textbf{i}_k}))$. By the Lemma \ref{wazny lemat} the jump of the deformation $(f_s^{\textbf{i}_k})$ is given by the formula
$$\lambda((f_s^{\textbf{i}_k}))=\lambda_k(p_k-1).$$
This formula implies that there is no "better" point in the hyperplane\\ $Ox_1\ldots x_{k-1}x_{k+1}\ldots x_n$ which could give a monomial deformation of $f_0$ with smaller jump. In fact, for any point $\textbf{i}\in J\cap Ox_1\ldots x_{k-1}x_{k+1}\ldots x_n$ we would have
$$\lambda((^kf_s^{\textbf{i}}))\geqslant \lambda_k$$
and hence again by the Lemma \ref{wazny lemat}
$$\lambda((f_s^{\textbf{i}}))=\lambda((^kf_s^{\textbf{i}}))(p_k-1)\geqslant \lambda_k(p_k-1).$$

Hence doing the same for all hyperplanes $Ox_1\ldots x_{k-1}x_{k+1}\ldots x_n$, $k=1,\ldots,n$, we get that the minimum of the jumps of monomial deformations associated to points lying in these hyperplanes is equal to
$$\min_{1\leqslant k\leqslant n}\lambda_k(p_k-1).$$
Denote this number by $\lambda^{\operatorname{hyp}}$.

\medskip
B. Now we consider monomial deformations of $f_0$ associated to points $\textbf{i}\not=0$ under $\Gamma_+(f_0)$ not lying in coordinate hyperplanes i.e. lying in the interior of the simplex $\Delta(P_0,P_1,\ldots,P_n)$, where $P_0=(0,\ldots,0)$, $P_1=(p_1,0,\ldots,0)$,$\ldots$, $P_n=(0,\ldots,0,p_n)$. Such a point $\textbf{i}=(i_1,\ldots,i_n)$ satisfies the conditions:
\begin{enumerate}[label=\arabic*.]
	\item\label{warunek 1} $0<i_k<p_k$, $k=1,\ldots n$,
	\item\label{warunek 2} $\dfrac{i_1}{p_1}+\ldots+\dfrac{i_n}{p_n}<1$ (or equivalently $i_1p_2\ldots p_n+\ldots+i_np_1\ldots p_{n-1} <p_1\ldots p_n$).
\end{enumerate}
For such a point $\textbf{i}$ we obviously have
$$\lambda((f_s^{\textbf{i}}))=n!\operatorname{vol}(\Delta(\textbf{i},P_1,\ldots,P_n))=p_1\ldots p_n-i_1p_2\ldots p_n-\ldots -i_np_1\ldots p_{n-1}.$$
This follows from the fact that $\operatorname{vol}(\Delta(P_0,\ldots,P_n))=\dfrac{p_1\ldots p_n}{n!}$ and \\ $\operatorname{vol}(\Delta(P_0,\ldots,P_{k-1},\textbf{i},P_{k+1},\ldots,P_n))=\dfrac{i_kp_1\ldots p_{k-1}p_{k+1}\ldots p_n}{n!}$, $k=1,\ldots,n$. So the problem is to find minimal $l$ in the set $\{1,\ldots,\lambda^{\operatorname{hyp}}-1\}$ for which there exist solutions $i_1,\ldots,i_n$ of the diophantine equation
\begin{equation}\label{rownanie diofantyczne z problemu}
i_1p_2\ldots p_n+\ldots+i_np_1\ldots p_{n-1}=p_1\ldots p_n -l
\end{equation}
satisfying conditions \ref{warunek 1} and \ref{warunek 2} First, we transform (\ref{rownanie diofantyczne z problemu}) into the more useful form
\begin{equation}
-i_1p_2\ldots p_n-\ldots-i_{n-1}p_1\ldots p_{n-2}p_n+(p_n-i_n)p_1\ldots p_{n-1}=l.
\end{equation}
If we denote $p'_k\assign p_1\ldots p_{k-1}p_{k+1}\ldots p_n$, $k=1,\ldots, n$, then we reduce the problem to find integer solutions 
$i_1,\ldots,i_{n-1},\widetilde{i_n}$ of the equations
\begin{equation}\label{ostateczna postac rownania z problemu}
-i_1p'_1-\ldots-i_{n-1}p'_{n-1}+\widetilde{i_n}p'_n=l,\quad l=1\ldots,\lambda^{\operatorname{hyp}}-1,
\end{equation}
satisfying
\begin{equation}\label{warunki do rownania z problemu}
0<i_k<p_k, \quad k=1,\ldots,n-1,\quad 0<\widetilde{i_n}<p_n.
\end{equation}
Then for such a minimal $l_0\in\{1,\ldots,\lambda^{\operatorname{hyp}}-1\}$ for which the solution of (\ref{ostateczna postac rownania z problemu}) 
with conditions 
exists we have
$$\lambda^{\operatorname{nd}}(f_0)=l_0.$$

We give a simple algorithm (using only Euclid's algorithm) to solve the equation (\ref{ostateczna postac rownania z problemu}) with conditions (\ref{warunki do rownania z problemu}) only in the case $\operatorname{GCD}(p_i,p_j)=1$, $i,j=1,\ldots,n$, $i\not=j$. In general case we also may compute $\lambda^{\operatorname{nd}}(f_0)$ but modulo finding solutions of (\ref{ostateczna postac rownania z problemu}) satisfying (\ref{warunki do rownania z problemu}).

Assume $f_0=z_1^{p_1}+\ldots+z_n^{p_n}$, $p_i\geqslant 2$, $\operatorname{GCD}(p_i,p_j)=1$, for $i\not =j$. First we consider the equation (\ref{ostateczna postac rownania z problemu}) for $l=1$. Since $p_1,\ldots,p_n$ are pairwise prime therefore $\operatorname{GCD}(p'_1,\ldots,p'_n)=1$. Hence by Euclid's algorithm we find integers $j_1,\ldots,j_n\in\mathbb{Z}$ such that
$$j_1p'_1+\ldots+j_np'_n=1.$$
By simple transformations we get another identity
\begin{equation}\label{rownanie difantyczne dla p_1,..,p_n wzglednie pierwszych}
-i_1p'_1-\ldots-i_{n-1}p'_{n-1}+i'_np'_n=1
\end{equation}
where $0<i_k<p_k$, $k=1,\ldots,n-1$, $i'_n\in\mathbb{N}$.
If $i'_n<p_n$ then for $\textbf{i}=(i_1,\ldots,i_{n-1},p_n-i'_n)$ we have $$\lambda((f^{\textbf{i}}_s))=1$$
and consequently
$$\lambda^{\operatorname{nd}}(f_0)=1.$$
If $i'_n\geqslant p_n$ then we claim there is no solution to our problem for $l=1$. Assume to the contrary that $j_1,\ldots,j_n$ is such a solution i.e.
$$-j_1p'_1-\ldots-j_{n-1}p'_{n-1}+(p_n-j_n)p'_n=1,\quad 0<j_k<p_k,\quad k=1,\ldots,n.$$
Hence and by (\ref{rownanie difantyczne dla p_1,..,p_n wzglednie pierwszych}) we obtain
$$(j_1-i_1)p'_1+\ldots+(j_{n-1}-i_{n-1})p'_{n-1}+(j_n+i'_n-p_n)p'_n=0.$$
Since $p_1|p'_k$, $k=2,\ldots,n$, then $p_1|(j_1-i_1)p'_1$. But $p_1$ is relatively prime to $p_2,\ldots,p_n$. Hence $p_1|(j_1-i_1)$. But $|j_1-i_1|<p_1$ which implies $j_1=i_1$. 
\newline Similarly $j_2=i_2,\ldots,j_{n-1}=i_{n-1}$. This implies $i'_n=p_n-j_n<p_n$, which contradicts our supposition $i'_n\geqslant p_n$. 

So in this case $\lambda^{\operatorname{nd}}(f_0)>1$ and we repeat our considerations for $l=2$ etc. up to $\lambda^{\operatorname{hyp}}-1$. If in one of this steps we find solutions $i_1,\ldots,i_n$ of (\ref{ostateczna postac rownania z problemu}) satisfying conditions (\ref{warunki do rownania z problemu}) we stop the procedure and in this case $\lambda^{\operatorname{nd}}(f_0)=l$ and the monomial deformation $(f^{\textbf{i}}_s)$ for $\textbf{i}=(i_1,\ldots,i_{n-1},p_n-i_n)$ realizes the jump. If the above search fails we conclude
$$\lambda^{\operatorname{nd}}(f_0)=\lambda^{\operatorname{hyp}}.$$
We may sum up the above considerations in the following theorem.
\begin{thm}
	Let $f_0 \in \mathcal{O}^{n}$ be a convenient and non-degenerate singularity with only one\\ $(n-1)$-dimensional face of its Newton diagram. Assume
	the vertices (on axes) of this face are $(p_1, 0,\ldots,0),\ldots,(0, \ldots, p_n), \;p_i\geqslant 2$. Then the non-degenerate jump $\uplambda^{\operatorname{nd}}(f_0)$ of $f_0$ can be inductively (with respect to $n$) calculated as follows. For $n=1$, $\uplambda^{\operatorname{nd}}(f_0)=1$. For $n\geqslant 2$ let $\uplambda_k$ be the non-degenerate jump of the restriction $^kf_0$ of $f_0$ to the hyperplane $\{z_k=0\}$ and $\uplambda^{\operatorname{hyp}}\assign \min\limits_{1\leqslant k\leqslant n}(\uplambda_k(p_k-1))$. Then putting $p'_k\assign p_1\ldots p_{k-1}p_{k+1}\ldots p_n$, $k=1,\ldots,n$ we have
	
	\begin{equation}\label{r:glowna formula}
	 \uplambda^{\operatorname{nd}} (f_0) = \left\{ \begin{array}{ll}
	\min\limits_{1\leqslant l<\uplambda^{\operatorname{hyp}}} (l) 
	
	& \text{if there exist } j_1,\ldots j_n\in\mathbb{Z} \text{ such that }\\
	& j_1p'_1+\ldots+j_np'_n=l,\; 0<-j_k<p_k,\vspace{2mm}\\
	
	&  k=1,\ldots,n-1,\; 0<j_n<p_n, \\
	
&\\
	\uplambda^{\operatorname{hyp}}  &\text{otherwise}.
	\end{array} \right. 
	\end{equation}
	In the first case the jump is realized by the monomial deformation $(f^{\textup{\textbf{i}}}_s)$ where $\textup{\textbf{i}}=(-j_1,\ldots,-j_{n-1},p_n-j_n)$ and in the second case by $(f^{\textup{\textbf{i}}}_s)$ where $\textup{\textbf{i}}$ is the point in a coordinate hyperplane $\{x_{k_0}=0\}$ for which $\uplambda^{\operatorname{hyp}}=\uplambda_{k_0}(p_{k_0}-1)$, realizing the jump for the restriction $^{k_0}f_0$ of $f_0$ to the hyperplane $\{z_{k_0}=0\}$.
	
	Moreover in the case $\operatorname{GCD}(p_i,p_j)=1$ for $i\not=j$ the jump $\uplambda^{\operatorname{nd}}(f_0)$ can be computed using only Euclid's algorithm. 
	
\end{thm}

\begin{ex}[Bodin \cite{Bod07}, Walewska \cite{Wal10}]\label{ex1} For $n=2$ and $f_0(\textbf{z})=z_1^{p_1}+z_2^{p_2}$ putting $d\assign\operatorname{GCD}(p_1,p_2)$ we have from the above theorem 
	
	\[ \uplambda^{\operatorname{nd}} (f_0) = \left\{ \begin{array}{lll}
	d &\text{if}& d<\min(p_1,p_2) \\
	d-1 &\text{if} & d=\min(p_1,p_2)
	\end{array}. \right. \]
	In fact, we have $\uplambda_1=\uplambda_2=1$, $\uplambda^{\operatorname{hyp}}=\min(p_1-1,p_2-1)$ and we always have the solution $j_1,j_2$ of $j_1p_2+j_2p_1=d$, $0<-j_1<p_1$, $0<j_2<p_2$.
	
\end{ex}

\begin{ex}
	Let $f_0(\textbf{z})=z_1^{11}+z_2^{6}+z_3^{5}$. We have $p_1=11$, $p_2=6$, $p_3=5$ and by  \textup{Example \ref{ex1}} $\uplambda_1=\uplambda_2=\uplambda_3=1$, $\uplambda^{\operatorname{hyp}}=4$. By the algorithm given in the proof we easily find solutions of diophantine equations in the formula (\ref{r:glowna formula}) for $l=1,2,3$
	$$-4p_2p_3-5p_1p_3+5p_1p_2=1,$$
	$$-8p_2p_3-4p_1p_3+7p_1p_2=2,$$
	$$-1p_2p_3-3p_1p_3+3p_1p_2=3.$$
	Only the last equation fulfills the condition in the formula. Hence
	$$\uplambda^{\operatorname{nd}}(f_0)=3$$
	and the jump is realized by the monomial deformation
	$$f_s^{(1,3,2)}(\textbf{z})=z_1^{11}+z_2^{6}+z_3^{5}+sz_1z_2^{3}z_3^{2}.$$
\end{ex}

\section*{Acknowledgement}
The first author was partially supported by the Polish
National Science Centre (NCN), Grant No. 2012/07/B/ST1/03293.

\bibliographystyle{plain}
\bibliography{bibliografia}

\end{document}